\documentclass[letterpaper, 10pt, conference]{ieeeconf}      % Use this line for a4
\IEEEoverridecommandlockouts                              % This command is only
\newcommand{\ud}{\,\mathrm{d}}

\usepackage{graphics} % for pdf, bitmapped graphics files
\usepackage{epsfig} % for postscript graphics files
\usepackage{mathptmx} % assumes new font selection scheme installed
\usepackage{times} % assumes new font selection scheme installed
\usepackage{amsmath} % assumes amsmath package installed
\usepackage{amssymb}  % assumes amsmath package installed
\usepackage{mathtools}
\usepackage{tabularx}
\usepackage{multirow}
\usepackage{csquotes}
\usepackage{color}
\usepackage{floatrow}
\usepackage{caption}
\usepackage{subfigure}  % For side-by-side figures

\usepackage{ifpdf}
\ifpdf
    \usepackage{graphicx}
    \usepackage{epstopdf}
\else
    \usepackage{graphicx}
\fi

\newtheorem{theorem}{Theorem}

\newtheorem{remark}{Remark}

\newtheorem{assumption}{Assumption}

\usepackage{algorithm}
\usepackage{algorithmic}

\newlength{\noteWidth}
\long\def\notes#1{\ifinner
             {\tiny #1}
             \else
              \marginpar{\parbox[t]{\noteWidth}{\raggedright\tiny #1}}
               \fi}

\def\notes#1{\typeout{#1 !!!}}  % for final copy

\makeatletter
\newcommand{\rom}[1]{\romannumeral #1}
\newcommand{\Rom}[1]{\expandafter\@slowromancap\romannumeral #1@}
\makeatother

\newcounter{rmnum}

\newcounter{anum}

\def\Sec#1{Sec.~\ref{#1}}

\def\IEEEQEDclosed{\mbox{\rule[0pt]{1.3ex}{1.3ex}}}

\def\qed{\nobreak\hfill\IEEEQEDclosed}
\def\UZ{\underline{\mathcal{Z}}}

\DeclareMathOperator{\Tr}{Tr}

\def\Re{\mathbb{R}}

\def\clV{{\cal V}}
\def\clW{{\cal W}}

\def\Cinf{{C^{\infty}}}

\newcommand{\lr}[2]{\langle #1, #2 \rangle}
\newcommand{\sk}[1]{[\,#1\,]_{\times}}
\newcommand{\degg}[1]{{#1}^{\circ}}

\def\Expect{{\sf E}}
\def\K{{\sf K}}
\def\x{{R}}

\def\k{{\sf k}}
\def\id{{I}}
\def\innov{\mathrm{I}}
\def\lap{\Delta_{\rho}}

\def\UZ{\mathcal{Z}_t}
\def\Dt{\Delta t}

\def\ge{g^{(\epsilon)}}

\def\ken{k^{(\epsilon,N)}}

\def\grad{\text{grad}}
\def\div{\text{div}}

\def\det{\text{det}}
\def\exp{\text{exp}}

\title{\LARGE \bf
Attitude Estimation with Feedback Particle Filter
}

\author{Chi Zhang, Amirhossein Taghvaei and Prashant G. Mehta %and Sean P. Meyn% <-this % stops a space
\thanks{Financial support from the NSF CMMI grants 1334987 and 1462773 is gratefully acknowledged.}% <-this % stops a space
\thanks{C.~Zhang, A.~Taghvaei and P. G.~Mehta are with the Coordinated
  Science Laboratory and the Department of Mechanical Science and
  Engineering at the University of Illinois at Urbana-Champaign (UIUC)
{\tt\scriptsize \{czhang54; taghvae2; mehtapg\}@illinois.edu}}%; mehtapg@illinois.edu}}% <-this
}

\begin{document}

\maketitle

\begin{abstract}
This paper presents theory, application, and comparisons of the feedback particle filter (FPF) algorithm for the problem 
of attitude estimation. The paper builds upon our recent work on the exact FPF solution of the 
continuous-time nonlinear filtering problem on compact Lie groups. 
% The FPF on Lie groups possesses a gain-feedback structure 
% and provides a coordinate-free description that satisfies the geometric constraints of the manifold. 
In this paper, the details of the FPF algorithm are presented for the problem 
of attitude estimation\,--\,a nonlinear filtering problem on $SO(3)$.
The quaternions are employed for computational purposes. The algorithm requires a numerical solution of the 
filter gain function, and two methods are applied for this purpose. 
Comparisons are also provided between the FPF and some popular algorithms for attitude estimation 
on $SO(3)$, including the invariant EKF, the multiplicative EKF, and the unscented Kalman filter. 
Simulation results are presented that help illustrate the comparisons.

\end{abstract}

\section{Introduction}
\label{sec:intro}

Attitude estimation is important to numerous fields including 
localization of mobile robots \cite{barrau2015cdc, barczyk2013invariant, hua2014implementation}, 
visual tracking of objects \cite{choi2011robust, kwon2007particle}, 
and navigation of spacecrafts \cite{crassidis2003unscented, carmi2009adaptive}. 
The mathematical problem of attitude estimation is a nonlinear filtering problem 
on a matrix Lie group, in particular the special orthogonal group $SO(3)$.
The design of attitude filters thus requires consideration of the geometry of the manifold. 

A number of attitude filters have been proposed and applied for the aerospace applications.
A majority of these filters are based on the extended Kalman filter (EKF), e.g. 
the additive EKF \cite{bar1985attitude, choukroun2006novel} and the 
multiplicative EKF \cite{lefferts1982kalman, markley2003attitude}. 
The EKF-based filters require a linearized model of the estimation error. 
Such a model is typically derived using one of the many 
three-dimensional attitude representations, e.g. the Euler angle \cite{itzhack87}, 
the rotation vector \cite{pittelkau2003rotation}, 
and the modified Rodrigues parameter \cite{crassidis1996attitude}. 
These representations have also been employed in the construction of 
unscented Kalman filters \cite{crassidis2003unscented, cheon2007unscented}. 
More recently, group-theoretic methods for attitude estimation have been explored. 
Deterministic nonlinear observers that respect the intrinsic geometry of the Lie groups have appeared in   
\cite{mahony2008TAC, khosravian2016automatica, wu2015hybrid, batista2014attitude, bohn2014almost}. 
A class of symmetry-preserving observers have been proposed to 
exploit certain invariance properties 
\cite{bonnabel2008symmetry, bonnabel2009observer}, leading to 
the invariant EKF algorithm \cite{bonnabel2009IEKF, barrau2015TAC, barczyk2015invariant}, 
the invariant ensemble EKF \cite{barrau2015TAC}, 
and the invariant particle filter \cite{barrau2014cdc} within the stochastic filtering framework. 
Filters based on certain variational formulations on Lie groups 
have also been investigated \cite{zamani2013tac, berger2015second, izadi2014rigid}.
Particle filters for attitude estimation include the 
bootstrap particle filter \cite{cheng2010particle, oshman2006attitude}, the marginalized particle filter 
\cite{wang2014steady}, and the Rao-Blackwellized particle filter \cite{vernaza2006rao}. 
For more comprehensive review and performance comparison of the various attitude filters, c.f., 
\cite{crassidis2007survey, zamani2013thesis, izadi2015comparison}. 
Some of these filters are also described in  
Sec. \ref{sec:other_filters} for the purpose of comparisons with the proposed FPF algorithm.

The feedback particle filter (FPF) is an exact algorithm for the solution of the continuous-time 
nonlinear filtering problem. 
The FPF algorithm was originally proposed in the Euclidean setting of $\Re^n$ \cite{Tao_TAC}. 
In a recent paper from our group, the FPF was extended to filtering on compact matrix Lie groups \cite{ACC2016}. 
The FPF is an intrinsic algorithm: 
The particle dynamics, expressed in their Stratonovich form, respect the geometric constraints 
of the manifold. 
The update step in FPF has a gain-feedback structure where the gain needs to be obtained numerically as a 
solution to a certain linear Poisson equation. When the gain function can be exactly computed, 
the FPF is an {\em exact} algorithm. 
In this case, in the limit of large number of particles, the empirical distribution of the particles 
exactly matches the posterior distribution of the hidden state.

\medskip

The contributions of this paper are as follows:

\smallskip

\noindent {$\bullet$ \bf FPF algorithm for attitude estimation.} 
The FPF algorithm is presented for the problem of attitude estimation. 
The explicit form of the filter is described with respect to both the rotation matrix and the quaternion coordinate, 
with the latter being demonstrated for computational purposes.

\smallskip

\noindent {$\bullet$ \bf Numerical solution of the gain function.} 
The FPF algorithm requires numerical approximation of the gain function as a solution 
to a linear Poisson equation on the Lie group. For this purpose, two numerical methods are proposed: 
In a Galerkin scheme, the gain function is approximated with a set of pre-defined basis functions. 
The second scheme involves solving a fixed-point equation associated with the weighted Laplacian 
operator on the manifold.

\smallskip

\noindent {$\bullet$ \bf Comparison of attitude filters.} 
For the purpose of comparison, the invariant EKF, the multiplicative EKF, and the UKF algorithms 
are briefly reviewed. Simulation studies are presented to compare performance between these 
filters and the proposed FPF algorithm.

\smallskip

The remainder of this paper is organized as follows: 
After a brief review in \Sec{sec:prelim} of the relevant Lie group preliminaries, 
the problem of attitude estimation is formulated in Sec. \ref{sec:problem}. 
The FPF algorithm on $SO(3)$ is described in \ref{sec:FPF_SO3},  
and some other attitude filters are briefly reviewed in Sec. \ref{sec:other_filters}. 
Numerical simulations are contained in Sec. \ref{sec:simulation_SO3}. 
% ll supplemental materials appear in the appendix.

\section{Mathematical Preliminaries}
\label{sec:prelim}

\noindent{\bf Geometry of $SO(3)$:}
The special orthogonal group $SO(3)$ is the group of $3\times3$ 
matrices $R$ such that $RR^T=\id$ and $\det(R)=1$. 
% Throughout the paper, $\id$ denotes the identity matrix of appropriate dimension.
The Lie algebra $so(3)$ is the 3-dimensional 
inner product space of skew-symmetric matrices. 
The inner product is denoted as $\lr{\cdot}{\cdot}_{so(3)}$.
Given an orthonormal basis $\{E_1,E_2, E_3\}$, 
a vector $\omega = (\omega_1,\omega_2,\omega_3) \in \Re^3$ is uniquely mapped to an element 
in $so(3)$, denoted as $\sk{\omega}:=\omega_1E_1 + \omega_2E_2 + \omega_3E_3$. 
The exponential map of $\Omega \in so(3)$ is denoted as $\exp(\Omega)$, 
and the space of smooth real-valued functions $f:SO(3)\rightarrow\Re$ is
denoted as $\Cinf(G)$, where we write $G$ interchangeably as $SO(3)$.

\smallskip

\noindent {\bf Vector field:}
The Lie algebra is identified with the tangent space at the
identity matrix $\id \in SO(3)$, and used to construct a basis $\{E^R_1,E^R_2, E^R_3\}$ 
for the tangent space at $\x \in SO(3)$, where $E^R_n := R E_n$ for $n=1,2,3$.
Therefore, a smooth vector field, denoted as $\clV$, is expressed as,
  $$ \clV(\x) = v_1(\x) \, E_1^{\x} + v_2(\x) \, E_2^{\x} + v_3(\x) \, E_3^{\x}, $$
with $v_n(\x) \in \Cinf(G)$. We write $\clV = \x V$,
where $V(\x) := v_1(\x)\,E_1 + v_2(\x)\,E_2 + v_3(\x)\,E_3$ is an element of $so(3)$. 
The functions $\big(v_1(\x), v_2(\x), v_3(\x)\big)$ are called {\em coordinates} of $\clV$. 
The inner product of two vector fields is
\begin{equation*}
 \lr{\clV}{\clW}(\x) := \lr{V}{W}_{so(3)}(\x) = \sum_{n=1}^3 v_n(\x)w_n(\x).
\end{equation*}

With a slight abuse of notation, the action of the vector 
field $\clV$ on $f\in\Cinf(G)$ is denoted as,
\begin{equation*}
 V\cdot f(\x) := \frac{\ud}{\ud t}\Big|_{t=0} f \big(\x \, \exp(tV(\x))\big).
 \label{eq:action}
\end{equation*}
A smooth function, denoted as $\div\clV$, is then defined as, 
\begin{equation*}
 \div \clV(R) = \sum_{n=1}^3 E_n\cdot v_n(R).
 \label{eq:divergence}
\end{equation*}
We also define the vector field $\grad(\phi)$ for $\phi\in\Cinf(G)$ as,
\begin{equation*}
 \grad(\phi)(\x) = \x \, \K(\x),
 \label{eq:gradient_prelim}
\end{equation*}
where $\K(\x)\in so(3)$, with coordinates $\big(\k_1(\x),~ \k_2(\x),~ \k_3(\x)\big)$  
$:= \big(E_1\cdot\phi(\x),~ E_2\cdot\phi(\x),~ E_3\cdot\phi(\x)\big)$.
% The vector field acts on a function $f\in\Cinf(G)$ as,
% \begin{equation}
%  \K\cdot f(\x) = \sum_{n=1}^3 E_n\cdot\phi(\x) \, E_n\cdot f(\x) = \lr{\grad(\phi)}{\grad(f)}(\x). 
%  \label{eq:grad_inner_product}
% \end{equation}

\medskip

Apart from $\Cinf(G)$, we also consider the following function spaces: 
For a probability measure $\pi$ on $G$, $L^2(G;\pi)$ denotes the Hilbert space of 
functions on $G$ that satisfy $\pi(|f|^2) < \infty$ (\,here $\pi(|f|^2) := \int_G |f|^2 \ud\pi$\,);
$H^1(G;\pi)$ denotes the Hilbert space of functions $f$ such that $f$ and 
$E_n\cdot f$ (defined in the weak sense) are all in $L^2(G;\pi)$.

\medskip

\noindent{\bf Quaternions:} 
Quaternions provide a computationally efficient coordinate representation for $SO(3)$. 
A unit quaternion has the general form
\begin{align*}
  q & = (q_0, ~q_1, ~q_2, ~q_3) \notag \\ 
  & =\Big( \cos(\frac{\theta}{2}), ~\sin(\frac{\theta}{2})\omega_1, ~\sin(\frac{\theta}{2})\omega_2, ~
  \sin(\frac{\theta}{2})\omega_3 \Big), \label{eq:quat_def}
\end{align*}
and represents rotation of angle $\theta$ about the axis defined 
by the unit vector $(\omega_1, \omega_2, \omega_3)$.
As with $SO(3)$, the space of quaternions admits a Lie group structure: 
The identity quaternion is $q_I=(1,0,0,0)$, 
the inverse of $q$ is $q^{-1}=(q_0,-q_1,-q_2,-q_3)$, 
and the multiplication is defined as,
\begin{equation*}
 p \otimes q = 
 \begin{bmatrix}
  p_0q_0 - p_V \cdot q_V \\
  p_0q_V + q_0p_V + p_V \times q_V
 \end{bmatrix},
 \label{eq:quaternion_multiplication}
\end{equation*}
where $p_V = (p_1,p_2,p_3)$, $q_V = (q_1,q_2,q_3)$, and $\cdot$ and $\times$ denote the 
dot product and the cross product of two vectors.

Given a unit quaternion $q$, the corresponding rotation matrix $R = R(q)\in SO(3)$ is calculated by, 
\begin{equation}
 R = 
 \begin{bmatrix}
  2q_0^2 + 2q_1^2 - 1 & 2(q_1q_2-q_0q_3) & 2(q_1q_3+q_0q_2) \\
  2(q_1q_2+q_0q_3) & 2q_0^2 + 2q_2^2 - 1 & 2(q_2q_3-q_0q_1) \\
  2(q_1q_3-q_0q_2) & 2(q_2q_3+q_0q_1) & 2q_0^2 + 2q_3^2 - 1
 \end{bmatrix}.
\label{eq:convert_qR}
\end{equation}

\medskip

For more comprehensive introduction of Lie groups and quaternions, we refer the reader to 
\cite{chirikjian2000, trawny2005indirect}.

\section{Attitude estimation problem statement}
\label{sec:problem}

\subsection{Process model}
\label{sec:process_model}

A kinematic model of rigid body is given by,
\begin{equation}
 \ud R_t = R_t\Omega_t \ud t + R_t \circ \sk{\ud B_t},
 \label{eq:kinematics}
\end{equation}
where $R_t \in SO(3)$ is the orientation of the rigid body at time $t$, 
expressed with respect to an inertial frame, 
$\Omega_t = \sk{\omega_t}$ represents the angular velocity expressed in 
the body frame, and $B_t$ is a standard Wiener process in $\Re^3$. 
Both $\Omega_t$ and $\sk{\ud B_t}$ are elements of $so(3)$. 
The $\circ$ before $\ud B_t$ indicates that the stochastic differential equation (sde) 
\eqref{eq:kinematics} is expressed in its Stratonovich form. 
% A similar kinematic model, expressed in the It\^{o} form, appears 
% in \cite{solo2015cdc}.

Using the quaternion coordinate, \eqref{eq:kinematics} is written as,
\begin{equation}
 \ud q_t = \frac{1}{2} q_t \otimes (\omega_t \ud t + \ud B_t), 
 \label{eq:kinematics_quat}
\end{equation}
where, by a slight abuse of notation, $\omega_t \in \Re^3$ is interpreted as a quaternion 
$(0,\omega_t)$, and $\ud B_t$ is interpreted similarly. The sde \eqref{eq:kinematics_quat} is 
also interpreted in the Stratonovich sense.

\subsection{Measurement model}
\label{sec:measurement_model}

\noindent {\bf Accelerometer:} 
In the absence of translational motion, the accelerometer is modeled as, 
\begin{equation}
 \ud Z_t^{g} = R_t^T r^g \ud t + \ud W_t^g,
 \label{eq:accelerometer}
\end{equation}
where $r^g \in \Re^3$ is the 
unit vector in the inertial frame aligned with the gravity, and $W_t^g$ is a standard Wiener process in $\Re^3$.

\noindent {\bf Magnetometer:}
The model of the magnetometer is of a similar form,
\begin{equation}
 \ud Z_t^{b} = R_t^T r^b \ud t + \ud W_t^b,
 \label{eq:magnetometer} 
\end{equation}
where $r^b \in \Re^3$ is the 
unit vector in the inertial frame aligned with the local magnetic field, 
and $W_t^b$ is a standard Wiener process in $\Re^3$.

\subsection{Nonlinear filtering problem on $SO(3)$}
\label{sec:filtering_problem}

In terms of the process and measurement models, the nonlinear filtering problem for 
attitude estimation is succinctly expressed as,
\begin{subequations}
 \begin{align}
  \ud R_t & = R_t \Omega_t \ud t + R_t \circ \sk{\ud B_t}, \label{eq:kinematics_1} \\
  \ud Z_t & = h(R_t) \ud t + \ud W_t, \label{eq:observation}
 \end{align}
\end{subequations}
where $\Omega_t = \sk{\omega_t}$ is the angular velocity, 
$h:SO(3) \rightarrow \Re^m$ is a given nonlinear function whose $j$-th coordinate 
is denoted as $h_j$ (i.e. $h=(h_1,h_2,...,h_m)$), and $W_t$ is a standard Wiener process in $\Re^m$. 
Note that \eqref{eq:observation} encapsulates the sensor models given in 
\eqref{eq:accelerometer} and \eqref{eq:magnetometer} with a single equation. 
For the purpose of this paper, it is not necessary to assume that the models are linear.
It is assumed that $B_t$ and $W_t$ are mutually independent, and independent of the initial 
condition $R_0$ which is drawn from a known initial distribution, denoted as $\pi_0^{*}$. 
% By scaling, it is also assumed without loss of generality that 
% the covariance matrices associated with $\{B_t\}$ and $\{W_t\}$ are identity matrices.

The objective of the attitude estimation problem, described by \eqref{eq:kinematics_1} 
and \eqref{eq:observation}, is to compute the conditional 
distribution of $R_t$ given the history of measurements (filtration) $\UZ=\sigma(Z_s:s\leq t)$. 
The conditional distribution, denoted as $\pi_t^*$, acts on a function $f\in\Cinf(G)$ 
according to,
\begin{equation*}
 \pi_t^*(f) := \Expect[f(R_t) | \UZ].
 \label{eq:true_posterior}
\end{equation*}

\smallskip

\begin{remark}
 There are a number of simplifying assumptions implicit in the model defined in 
 \eqref{eq:kinematics_1} and \eqref{eq:observation}. In practice, $\omega_t$ needs to be 
 estimated from noisy gyroscope measurements and there is translational motion as well. 
 This will require additional models which can be easily incorporated within the proposed 
 filtering framework. 
 
 The purpose of this paper is to elucidate the geometric aspects 
 of the FPF in the simplest possible setting of $SO(3)$. More practical FPF-based filters 
 that also incorporate models for translational motion, measurements of $\omega_t$ 
 from gyroscope, effects of translational motion on accelerometer, and effects of sensor bias 
 are subject of separate publication.
 \qed
 \label{remark:gyro}
\end{remark}

\section{Feedback Particle Filter on $SO(3)$}
\label{sec:FPF_SO3}

\subsection{FPF on $SO(3)$}
\label{eq:FPF}

The feedback particle filter is a controlled system 
with $N$ stochastic processes $\{R_t^i\}_{i=1}^N$ where $R_t^i \in SO(3)$ 
\footnote{Although the rotation matrix parameterization of $SO(3)$ is used, the filter is intrinsic. 
The FPF using the quaternion appears in Sec. \ref{sec:quaternion_rep}.}.
The conditional distribution of the particle $R_t^i$ given $\UZ$ is denoted by $\pi_t$, 
which acts on $f\in\Cinf(G)$ according to,
\begin{equation*}
 \pi_t(f) := \Expect[f(R_t^i) | \UZ]. 
 \label{eq:particle_posterior}
\end{equation*}
% The associated density function is denoted by $\rho$.

The dynamics of the $i$-th particle is defined by,
\begin{equation}
 ~~\ud R_t^i = \underbrace{R_t^i \, \Omega_t \ud t + R_t^i \circ \sk{\ud B_t^i}}_{\text{propagation}} 
              + \underbrace{R_t^i \, \sk{\K(R_t^i,t) \circ \ud \innov_t^i}}_{\text{measurement update}},
 \label{eq:particle_dyn}
\end{equation}
where $\{B_t^i\}_{i=1}^N$ are mutually independent standard Wiener processes in $\Re^3$, 
and $R_0^i$ is drawn from the initial distribution $\pi_0^{*}$.  
The $i$-th particle implements the Bayesian update step\,--\,to account for the conditioning due to 
the measurements\,--\,as gain $\K(R_t^i)$ times an error $\ud \innov_t^i$. 
The resulting control input to the $i$-th particle is an element of the Lie algebra $so(3)$.

The error $\ud \innov_t^i$ is a modified form of the innovation process: 
\begin{equation}
 \ud \innov_t^i = \ud Z_t - \frac{1}{2} \big( h(R_t^i) + \hat{h} \big) \ud t,
 \label{eq:innovation}
\end{equation}
where $\hat{h} := \pi_t(h)$. In a numerical implementation, we approximate 
$\hat{h} \approx \frac{1}{N}\sum_{i=1}^N h(R_t^i) =: \hat{h}^{(N)}$.

The gain function $\K$ is a $3 \times m$ matrix whose entries are obtained as follows: 
For $j = 1,2,...,m$, the $j$-th column of $\K$ is the coordinate of the vector field $\grad (\phi_j)$, 
where the function $\phi_j \in H^1(G;\pi)$ is a solution to the Poisson equation,
\begin{equation}
 \begin{aligned}
  & \pi_t \big( \lr{\grad(\phi_j)}{\grad(\psi)} \big) = \pi_t \big( (h_j-\hat{h}_j) \psi \big), \\
  & \pi_t (\phi_j) = 0 ~~~~(\text{normalization}), 
 \end{aligned}
 \label{eq:BVP}
\end{equation}
for all $\psi \in H^1(G;\pi)$.
This linear partial differential equation (pde) has to be solved for each $j=1,2,...,m$, and for each time $t \geq 0$. 
The existence-uniqueness of the solution of \eqref{eq:BVP} requires additional assumptions on 
$\pi_t$; c.f., \cite{Laugesen_2015}. 

\medskip

\begin{assumption}\
The distribution $\pi_t$ is absolutely 
continuous with respect to the uniform (Lebesgue) measure on $SO(3)$ with a positive density function $\rho$. 
\label{assumption:pi}
\qed
\end{assumption}

\medskip

Two numerical schemes for approximating the solution of \eqref{eq:BVP} appear in Sec. \ref{sec:galerkin_gain} and 
Sec. \ref{sec:kernel_gain}, respectively.

For the FPF \eqref{eq:particle_dyn}-\eqref{eq:BVP}, the following result is proved in \cite{ACC2016} 
that relates $\pi_t$ to $\pi_t^{*}$: 

\smallskip

\begin{theorem}
 Consider the particle system that evolves according to \eqref{eq:particle_dyn}, 
 where the gain function is obtained as solution to the Poisson equation \eqref{eq:BVP}, and 
 the error is defined as in \eqref{eq:innovation}. 
 Suppose that Assumption \ref{assumption:pi} holds. 
 Then assuming $\pi_0=\pi_0^{*}$, we have 
 $$ \pi_t(f) = \pi_t^{*}(f),$$
 for all $t > 0$ and all function $f\in\Cinf(G)$.
\label{thm:consistency} 
\qed
\end{theorem}

\subsection{Quaternion representation}
\label{sec:quaternion_rep}

For numerical purposes, it is convenient to express the FPF with respect to 
the quaternion coordinate. In this coordinate, the dynamics of the $i$-th particle evolves according to,
\begin{equation}
 \ud q_t^i = \frac{1}{2} \, q_t^i \otimes \ud \nu_t^i, 
 \label{eq:FPF_quat}
\end{equation}
where $q_t^i$ is the quaternion state of the $i$-th particle, and $\nu_t^i \in \Re^3$ evolves according to,
\begin{equation}
 ~~\ud\nu_t^i = \omega_t \ud t + \ud B_t^i + \K(q_t^i) 
 \circ \Big( \ud Z_t-\frac{h(q_t^i) + \hat{h}}{2}\, \ud t \Big),
 \label{eq:FPF_domega}
\end{equation}
where $\K(q,t)=\K(R(q),t)$ and $h(q)=h(R(q))$, with $R=R(q)$ given by the formula \eqref{eq:convert_qR}.

\subsection{Galerkin gain function approximation}
\label{sec:galerkin_gain}

In this section, a Galerkin scheme is presented to approximate the solution of the Poisson equation \eqref{eq:BVP}. 
Since the equations for each $j=1,2,...,m$ are uncoupled, without loss of generality, 
a scalar-valued measurement is assumed (i.e., $m=1$, and $\phi_j$, $h_j$ are 
denoted as $\phi$, $h$). As the time $t$ is fixed, 
the explicit dependence on $t$ is suppressed (i.e., we denote $\pi_t$ as $\pi$, $R_t^i$ as $R^i$). 
This notation is also used in Sec. \ref{sec:kernel_gain}.

In a Galerkin scheme, the solution $\phi$ is approximated as,
 $$ \phi = \sum_{l=1}^L \kappa_l \, \psi_l, $$
where $\{\psi_l\}_{l=1}^L$ is a given (assumed) set of {\em basis functions} on $SO(3)$. 
The gain function $\K = (\k_1,\k_2,\k_3)$, defined as the coordinates of $\grad(\phi)$, is then given by,
\begin{equation*}
 \k_n = \sum_{l=1}^L \kappa_l \, E_n\cdot \psi_l,~~~n = 1,2,3.
 \label{eq:gain}
\end{equation*}

The finite-dimensional approximation of the Poisson equation \eqref{eq:BVP} 
is to choose coefficients $\{\kappa_l\}_{l=1}^L$ such that,
\begin{equation}
  \sum_{l=1}^L \kappa_l \, \pi\big( \lr{\grad(\psi_l)}{\grad(\psi)} \big)
     = \pi \big( (h-\hat{h})\psi \big),
  \label{eq:FEM}
\end{equation}
for all $\psi \in \text{span}\{\psi_1,...,\psi_L\} \subset H^1(G;\pi)$. 
On taking $\psi=\psi_1,...,\psi_L$, \eqref{eq:FEM} 
is compactly written as a linear matrix equation,
\begin{equation}
 A\kappa = b,
 \label{eq:kappa}
\end{equation}
where $\kappa := (\kappa_1,\hdots,\kappa_L)$. 
The $L\times L$ matrix $A$ and the $L\times 1$ vector $b$ are defined and approximated as,
\begin{align}
 [A]_{kl} & = \pi\big(\lr{\grad(\psi_l)}{\grad(\psi_k)} \big) \notag \\
 & \approx \frac{1}{N} \sum_{i=1}^N \lr{\grad(\psi_l)({\x}^i)}{\grad(\psi_k)({\x}^i)} \notag \\
 & = \frac{1}{N} \sum_{i=1}^N \sum_{n=1}^3 (E_n\cdot\psi_l)({\x}^i) \, (E_n\cdot\psi_k)({\x}^i), \label{eq:A} \\
 b_k & = \pi \big((h-\hat{h})\psi_k \big) 
 \approx \frac{1}{N}\sum_{i=1}^N (h({\x}^i)-\hat{h})\psi_k({\x}^i), \label{eq:b}
\end{align}
where recall $\hat{h}\approx\frac{1}{N}\sum_{i=1}^N h({\x}^i) =: \hat{h}^{(N)}$.

Note that both the Poisson equation \eqref{eq:BVP} as well as its Galerkin finite-dimensional approximation 
\eqref{eq:kappa} are coordinate-free representations. Particle-based approximation of 
\eqref{eq:kappa}, viz. \eqref{eq:A} and \eqref{eq:b}, may be obtained using $R$ or $q$, 
or any other coordinate representation.

The choice of basis function is crucial in the Galerkin scheme, and one choice appears  
in Appendix \ref{sec:basis_fn}.

\subsection{Kernel-based gain function approximation}
\label{sec:kernel_gain}

In a kernel-based scheme, the solution to the Poisson equation \eqref{eq:BVP} 
is the solution of the following fixed-point equation for fixed positive $\tau$, 
\begin{equation}
 \phi = e^{\,\tau \, \lap} \phi + \int_{0}^{\tau} e^{\,s \, \lap} (h-\hat{h}) \ud s,
 \label{eq:fixed_point}
\end{equation}
where $e^{\,\tau\,\lap}$ is the semigroup associated with the {\em weighted Laplacian} on $SO(3)$, 
defined as $\lap:= (1/\rho)\,\div\big(\rho\,\grad(\phi)\big)$, 
where $\rho$ is the density of $\pi$.  
For small time $\tau=\epsilon$, 
the operator $e^{\,\tau\,\lap}$ is approximated using the particles as,  
\begin{equation}
 e^{\epsilon\,\lap} \phi(R) \approx \frac{\frac{1}{N}\sum_{i=1}^N \ken(R, R^i) \phi(R^i)}
 {\frac{1}{N}\sum_{i=1}^N \ken(R, R^i)}, 
 \label{eq:semigroup_approx}
\end{equation}
where the kernel $\ken: SO(3) \times SO(3) \rightarrow \Re$ is given by,
\begin{equation}
 \ken (R_1, R_2) = \frac{\ge(R_1, R_2)}
 {\sqrt{\frac{1}{N}\sum_{i=1}^N \ge(R_1, R^i)} 
 \sqrt{\frac{1}{N}\sum_{i=1}^N \ge(R_2, R^i)}}, 
 \label{eq:modified_kernel}
\end{equation}
and the Gaussian kernel $\ge$ is defined as, 
\begin{equation}
 \ge (R_1, R_2) := \frac{1}{(4\pi\epsilon)^{3/2}} \exp \Big( -\frac{|R_1-R_2|_F^2}{4\epsilon} \Big),
 \label{eq:kernel} 
\end{equation}
where $\epsilon$ is a small positive parameter, and $|\cdot|_F$ denotes the Frobenius norm of a matrix. 
The justification for the approximation \eqref{eq:semigroup_approx} appears in \cite{coifman2006diffusion}.

The approximation \eqref{eq:semigroup_approx} yields a finite-dimensional approximation of 
the fixed-point equation \eqref{eq:fixed_point}: 
\begin{equation}
 \varPhi = T^{(N)} \varPhi + \epsilon H^{(N)},
 \label{eq:fixed_point_approx}
\end{equation}
where $\varPhi \in \Re^N$ is the approximate solution that needs to be computed, 
$H^{(N)} = \big( h(R^1)-\hat{h}^{(N)}, h(R^2)-\hat{h}^{(N)}, ...,h(R^N)-\hat{h}^{(N)} \big)$, 
and $T^{(N)}\in\Re^{N\times N}$ whose entries are given by,
\begin{equation}
 T^{(N)}_{ij} = \frac{\ken (R^i, R^j)}{\sum_{l=1}^N \ken(R^i,R^l)}.
 \label{eq:Tij}
\end{equation}
Note that $T^{(N)}$ is a stochastic matrix with positive entries, and as a result, 
the fixed-point equation \eqref{eq:fixed_point_approx} is a contraction 
on the space of normalized vectors 
% (i.e., $\sum_{k=1}^N \varPhi_k = 0$\,). 
The solution can be obtained by successive approximations. 
% The number of iterations can be reduced is the 
% initial condition equals the solution at the previous time step. 
The solution $\phi$ of \eqref{eq:fixed_point}, evaluated at the particles, 
is then approximated as $\phi(R^i) \approx \varPhi_i$, the $i$-th entry of $\varPhi$. 

The gain function is given by $\K=(\k_1,\k_2,\k_3)$, where $\k_n = E_n\cdot\phi$ for $n=1,2,3$, 
and is evaluated at the particles according to,
% \begin{equation}
%  \begin{aligned}
%   E_n\cdot\phi(R^i) = & -\epsilon \, E_n R^{i\,T} r \\
%   & + \frac{1}{2\epsilon} \Big[ \big( S_n\varPhi \big)_i - 
%   \big( S_n \mathbf{1} \big)_i\big( T^{(N)}\varPhi \big)_i \Big]
%  \end{aligned}
%  \label{eq:kernel_grad}
% \end{equation}
\begin{equation}
 E_n\cdot\phi(R^i) = -\epsilon \, E_n \cdot h(R^i) 
 + \frac{1}{2\epsilon} \big[ \big( S_n\varPhi \big)_i - 
  \big( S_n \mathbf{1} \big)_i\big( T^{(N)}\varPhi \big)_i \big],
 \label{eq:kernel_grad}
\end{equation}
where $\mathbf{1} = (1,1,...,1)\in\Re^N$, and the entries of the $N\times N$ matrix $S_n$ are given by, 
  $$ (S_n)_{ij} = T^{(N)}_{ij}\,\Tr(R^iE_nR^j), $$
where $\Tr(\cdot)$ denotes the trace of a matrix. 

\medskip

\begin{remark}
 The theory for the kernel-based gain function approximation, together with its 
 convergence analysis and numerical illustration, appears in a companion paper \cite{AmirCDC2016}.
 \qed
\end{remark}

\subsection{FPF algorithm}
\label{sec:FPF_algorithm}

The FPF algorithm is numerically implemented using the quaternion coordinate, 
and is described in Algorithm 1. 
The algorithm simulates $N$ particles, $\{q_t^i\}_{i=1}^N$, according to the sde's \eqref{eq:FPF_quat} and \eqref{eq:FPF_domega}, 
with the initial conditions $\{q_0^i\}_{i=1}^N$ sampled i.i.d. from a given prior distribution $\pi_0^{*}$. 
The gain function is approximated using either the Galerkin scheme (see Sec. \ref{sec:galerkin_gain} 
and Algorithm 2), or the kernel-based scheme (see Sec. \ref{sec:kernel_gain} and Algorithm 3). 

Given a particle set $\{q_t^i\}_{i=1}^N$, its empirical mean is obtained as the 
eigenvector (with norm 1) of the $4\times4$ matrix $Q = \frac{1}{N} \sum\nolimits_{i=1}^N q_t^i {q_t^{i~T}}$, 
corresponding to its largest eigenvalue 
\cite{markley2007averaging}.

\begin{algorithm}[H]
        \label{alg:FPF-LG}
        \caption{Feedback Particle Filter on $SO(3)$}
        \begin{algorithmic}[1]
            \STATE {\bf initialization:}  sample $\{q_0^i\}_{i=1}^N$ from $\pi_0^{*}$
            \STATE Assign $t=0$
            \STATE {\bf iteration:} from $t$ to $t+\Dt$
            \STATE Calculate $\hat{h}^{(N)} = (1/N)\sum_{i=1}^N h(q_{t}^{i})$
	    \FOR{$i=1$ to $N$}
	      \STATE Generate a sample, $\Delta B_t^i$, from $N(0,\id)$
	      \STATE Calculate the error
		  \begin{equation*}
		      \Delta \innov_{t}^i := \Delta Z_{t} - (1/2)(h(q_t^i) + \hat{h}^{(N)}) \, \Dt
		      \label{eq:innovation_algo}
		  \end{equation*}
	      \STATE Calculate gain function $\K(q_t^i,t)$
	      \STATE Calculate $ \Delta \nu_t^i = \omega_t \, \Dt + \sqrt{\Dt} \, \Delta B_t^i + \K(q_t^i,t) \, \Delta \innov_{t}^i $
	      \STATE Propagate the particle $q_t^i$ according to
	             \begin{equation*}
		      q_{t+\Dt}^i = q_t^i \otimes 
		      \begin{bmatrix}
		      \cos \big( |\Delta\nu_t^i|/2 \big) \\
		      \frac{\Delta\nu_t^i}{|\Delta\nu_t^i|}\, \sin \big( |\Delta\nu_t^i|/2 \big)
		      \end{bmatrix}
		      \label{eq:dq_integration}
		     \end{equation*}
	             ($|\cdot|$ denotes the Euclidean norm in $\Re^3$)
	    \ENDFOR
	        \STATE Define matrix $Q = \frac{1}{N} \sum\nolimits_{i=1}^N q_{t+\Dt}^i {q_{t+\Dt}^{i~T}}$
            \STATE {\bf return:} empirical mean of $\{q_{t+\Dt}^i\}_{i=1}^N$, 
                                 i.e., the eigenvector of $Q$ associated with its largest eigenvalue 
            \STATE Assign $t = t+\Dt$
        \end{algorithmic}
\end{algorithm}

\begin{algorithm}[H]
        \label{alg:galerkin}
        \caption{Galerkin gain function approximation}
        \begin{algorithmic}[1]
        \STATE {\bf input:} Particles $\{q^i\}_{i=1}^N$
        \STATE Calculate $\hat{h}^{(N)} = (1/N)\sum_{i=1}^N h(q^{i})$
        \FOR{$k=1$ to $L$}
            \STATE Calculate $b_k = \frac{1}{N}\sum_{i=1}^N \big( h({q}^i)-\hat{h}^{(N)} \big)\psi_k({q}^i)$
            \FOR{$l=1$ to $L$}
                \STATE Calc. $A_{kl} = 
                \frac{1}{N} \sum_{i=1}^N \sum_{n=1}^3 (E_n\cdot\psi_l)({q}^i) \, (E_n\cdot\psi_k)({q}^i)$
            \ENDFOR
        \ENDFOR
        \STATE Solve the matrix equation $A\kappa=b$, with $A=[A_{kl}]$, $b=[b_k]$
        \STATE Calculate $\k_n(q^i) = \sum_{l=1}^L \kappa_l \, E_n\cdot\psi_l(q^i)$, for $n=1,2,3$
        \STATE {\bf return:} $\big\{\K(q^i) = \big( \k_1(q^i),\,\k_2(q^i),\,\k_3(q^i) \big)\big\}_{i=1}^N$
        \end{algorithmic}
\end{algorithm}

\section{Review of Some Attitude Filters}
\label{sec:other_filters}

In this section, we restrict our attention to the attitude estimation problem with linear observations of the 
form $h(R_t) = R_t^T r$ where $r$ is a known reference vector in the inertial frame 
(see the models of accelerometer and magnetometer in \eqref{eq:accelerometer}, \eqref{eq:magnetometer}). 
A majority of the literature deals with such linear models.
For discrete-time filters, it is convenient to define $Y_t := \frac{\ud Z_t}{\ud t}$, 
whose model is formally expressed as,
\begin{equation*}
 Y_t = R_t^T r + \dot{W}_t,
\end{equation*}
where $\dot{W}_t$ is a white noise process in $\Re^3$. 
In this section, we assume without loss of generality that the covariance matrix associated with $\dot{W}_t$ 
is the identity matrix.

The sequence of sampling instants is denoted as $\{t_n\}$, $n = 0,1,2,...$, 
with uniform time step $\Delta t = t_{n+1}-t_n$. 
The discrete-time sampled measurements are denoted as $\{Y_n\}$. Similarly, 
$\{R_n\}$ and $\{\omega_n\}$ denote the discrete-time samples of $R_t$ and $\omega_t$. 
Furthermore, $\widehat{R}_n$ denotes the posterior filter estimate at time $t_n$, 
$\widehat{R}_{n|n-1}$ denotes the filter estimate after the propagation step but 
before the measurement update, 
and $\Sigma_{n|n-1}$, $\Sigma_n$ denote the associated covariance matrices. 

\begin{algorithm}[t!]
        \label{alg:kernel}
        \caption{Kernel-based gain function approximation}
        \begin{algorithmic}[1]
        \STATE {\bf input:} Particles $\{q^i\}_{i=1}^N$, parameters $\epsilon$, $K$
        \STATE Calculate $\hat{h}^{(N)} = (1/N)\sum_{i=1}^N h(q^{i})$
        \FOR{$i=1$ to $N$}
            \STATE Calculate $H^{(N)}_i = h(q^i) - \hat{h}^{(N)}$
            \FOR{$j=1$ to $N$}
            		\STATE Calculate $\ge(q^i,q^j)$, $\ken(q^i,q^j)$ by \eqref{eq:kernel}, \eqref{eq:modified_kernel}
                \STATE Calculate $T^{(N)}_{ij}$ according to \eqref{eq:Tij}
                \FOR{$n=1,2,3$}
                    \STATE Calculate $(S_n)_{ij} = T^{(N)}_{ij} \Tr\big(R(q^i)\,E_n\,R(q^j)\big)$
                \ENDFOR
            \ENDFOR 
        \ENDFOR
        \STATE Assign $\varPhi^0$ as solution of \eqref{eq:fixed_point_approx} 
               in previous time step
        \FOR{$k=0$ to $K-1$}
	    \STATE Calculate $\varPhi^{k+1} = T^{(N)} \varPhi^k + \epsilon H^{(N)}$, with $T^{(N)}=[T^{(N)}_{ij}]$
	\ENDFOR
% 	\STATE Assign $\phi(q^i) = (\varPhi_K)_i$
	\FOR{$i=1$ to $N$}
	    \FOR{$n=1,2,3$}
% 		\STATE Calculate $\k_n(q^i)$ according to \eqref{eq:kernel_grad} with $\varPhi=\varPhi_K$
	        \STATE Calc. $\k_n(q^i)$ by \eqref{eq:kernel_grad} with $S_n=[(S_n)_{ij}]$ and $\varPhi=\varPhi^K$
	    \ENDFOR
	\ENDFOR
        \STATE {\bf return:} $\big\{\K(q^i) = \big( \k_1(q^i),\,\k_2(q^i),\,\k_3(q^i) \big)\big\}_{i=1}^N$
        \end{algorithmic}
\end{algorithm}

\subsection{Invariant extended Kalman filter}
\label{sec:IEKF}

The invariant EKF (IEKF) models the attitude at time $t_n$ as the product
\begin{equation}
 R_n = \delta R_n \, \widehat{R}_n,
 \label{eq:IEKF_prod}
\end{equation}
where the estimation error, $\delta R_n \in SO(3)$, is  
represented as $\delta R_n = \exp (\sk{\eta_n})$ where $\eta_n \in \Re^3$. 
% The error thus defined is shown to possess certain invariance property \cite{bonnabel2009IEKF}, 
% rendering the IEKF an intrinsic algorithm. 
At each time step, the estimate of $\eta_n$, denoted as $\hat{\eta}_n$, is obtained as follows: 

\noindent {\em (\rom{1}) Propagation step:}
\begin{align*}
 \widehat{R}_{n|n-1} & = \widehat{R}_{n-1} \exp(\sk{\omega_{n-1} \Delta t_{}}),  \\
 \Sigma_{n|n-1} & = \Sigma_{n-1} + (\Dt)\id. 
\end{align*}

\noindent {\em (\rom{2}) Update step:}
The innovation error is defined in the inertial frame,
\begin{equation*}
 \innov_n = \widehat{R}_{n|n-1} Y_n - r,
 \label{eq:IEKF_innovation}
\end{equation*}
and the gain matrix $\K_n$ is calculated according to,
\begin{equation*}
 \K_n = \Sigma_{n|n-1} H^T \big( H \Sigma_{n|n-1} H^T + \id \big)^{-1},
 \label{eq:IEKF_gain}
\end{equation*}
where $H = \sk{r}$. 

\noindent {\em (\rom{3}) Posterior update:}
\begin{align*}
 \hat{\eta}_n & = \K_n \innov_n,  \\
 \widehat{R}_n & = \exp\,(\sk{\hat{\eta}_n}) \, \widehat{R}_{n|n-1},  \\
 \Sigma_n & = (\id - \K_n H) \, \Sigma_{n|n-1}. 
\end{align*}

For more details of the IEKF algorithm, we refer the reader to \cite{barrau2015TAC}.

\subsection{Multiplicative extended Kalman filter}
\label{sec:MEKF}

The multiplicative EKF (MEKF) models the attitude at time $t_n$ as the product
\begin{equation}
 R_n = \widehat{R}_n \, \delta R_n,
 \label{eq:MEKF_prod}
\end{equation}
where the estimation error, $\delta R_n \in SO(3)$, 
is parameterized by some coordinate, 
e.g., the modified Rodrigues parameter \cite{crassidis1996attitude}. 
The coordinate is denoted as $a_n \in \Re^3$, and $\delta R_n = \delta R(a_n)$. 
At each time step, the estimate of $a_n$, denoted as $\hat{a}_n$, is obtained as follows:

\noindent {\em (\rom{1}) Propagation step:}
\begin{align*}
 \widehat{R}_{n|n-1} & = \widehat{R}_{n-1} \exp(\sk{\omega_{n-1} \Delta t_{}}),  \\
 \Sigma_{n|n-1} & = \Lambda \Sigma_{n-1} \Lambda^T + (\Delta t_{}) \id, 
\end{align*}
where $\Lambda = \id - \sk{\omega_{n-1}\Delta t}$. 

\noindent {\em (\rom{2}) Update step:}
\begin{equation*}
 \innov_n = Y_n - \widehat{R}_{n|n-1}^T r,
 \label{eq:MEKF_innovation}
\end{equation*}
\begin{equation*}
 \K_n = \Sigma_{n|n-1} H^T \big( H \Sigma_{n|n-1} H^T + \id \big)^{-1},
 \label{eq:MEKF_gain}
\end{equation*}
where $H = \sk{\widehat{R}_{n|n-1}^T r}$. In contrast to the IEKF, the 
innovation error in the MEKF is defined in the body frame.

\noindent {\em (\rom{3}) Posterior estimate:}
\begin{align*}
 \hat{a}_n & = \K_n \innov_n,  \\
 \widehat{R}_n & =\widehat{R}_{n|n-1} \, \delta R(\hat{a}_n),  \\
 \Sigma_n & = (\id - \K_n H) \, \Sigma_{n|n-1}. 
\end{align*}

For more details of the MEKF algorithm, we refer the reader to \cite{markley2003attitude, trawny2005indirect}.

\subsection{Unscented Kalman filter}
\label{sec:UKF}

The unscented Kalman filter (UKF) for attitude estimation, presented in \cite{crassidis2003unscented}, 
also uses the parameterization of the MEKF, i.e., 
  $$ R_n = \widehat{R}_n \delta R(a_n), $$
where $a_n \in \Re^3$ is the chosen coordinate. The estimate of $a_n$ is obtained by using a standard UKF 
in $\Re^3$. For equations of the algorithm, we refer the reader to \cite{crassidis2003unscented}. 

\subsection{Other filters}
 
Apart from the above, other types of attitude filters include the continuous-time 
IEKF \cite{bonnabel2009IEKF}, the geometric approximate minimum-energy (GAME) filter 
\cite{zamani2013tac}, and the bootstrap particle filter \cite{cheng2010particle}. 
These filters are not included in the simulation-based comparisons that are presented next.

\section{Simulations}
\label{sec:simulation_SO3}

\begin{figure*}
  \centering
%   \begin{tabular}{cc}
    \subfigure[Initial distribution: $\degg{30}$, sensor noise: $\degg{10}$]{
    \includegraphics[scale=0.4]{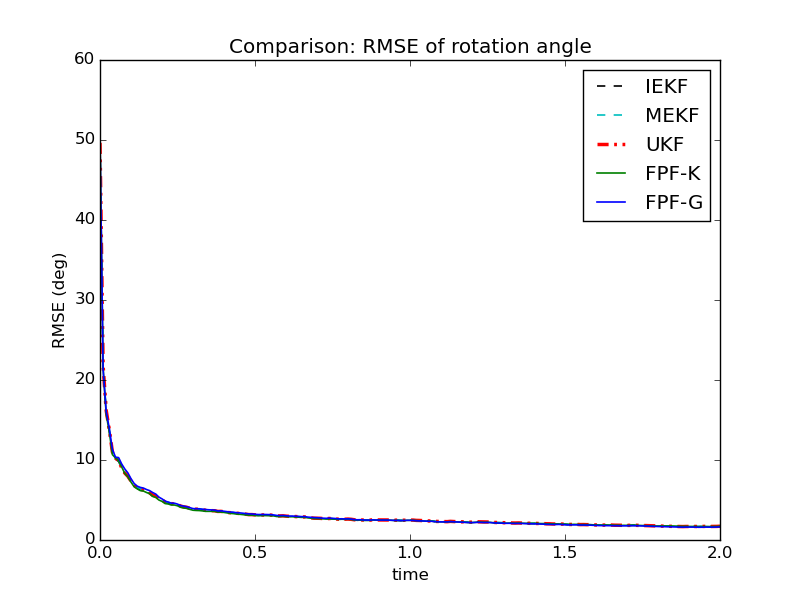}}
    \subfigure[Initial distribution: $\degg{60}$, sensor noise: $\degg{30}$]{
    \includegraphics[scale=0.4]{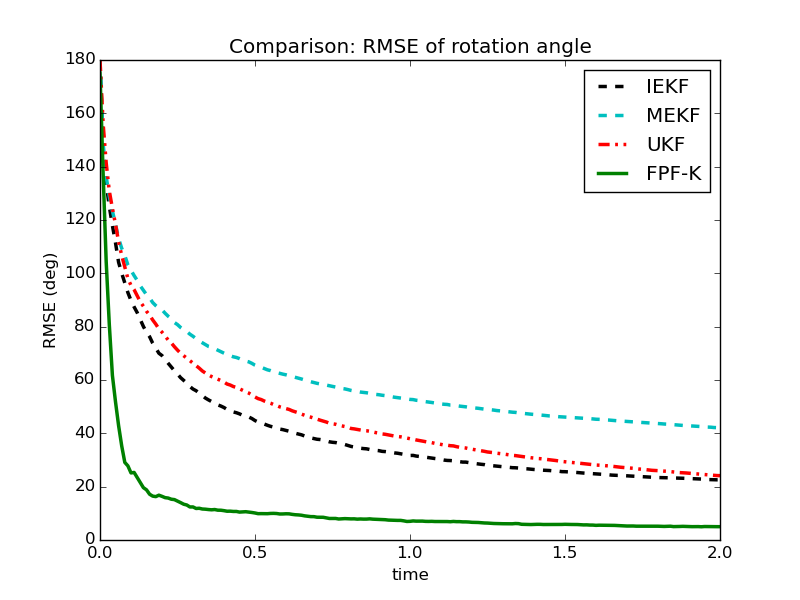}}
%     \end{tabular}
  \caption{Fig. 1: Comparison of filter performance with different initial distribution error and sensor noise.}
  \label{fig:init_var}
\end{figure*}

% \begin{figure*}
%   \centering
%     \subfigure[Initial dist. std.: $\degg{60}$, sensor noise std.: $\degg{30}$]{
%     \includegraphics[scale=0.38]{figures/ang_LA60_W30_N200_MC100_T200_discretization.png}}
%     \subfigure[Initial error: $\degg{180}$, sensor noise std.: $\degg{30}$]{
%     \includegraphics[scale=0.38]{figures/ang_A180_W30_N200_MC100_T200_discretization.png}}
%   \caption{Filter performance with large sensor noise and large initial estimation error.}
%   \label{fig:large}
% \end{figure*}

For numerical simulations of the filters, 
we consider the following attitude estimation problem,
\begin{align*}
  \ud q_t & = \frac{1}{2} q_t \otimes \big( \omega_t \ud t + \Sigma_B \, \ud B_t \big),  \\
  \ud Z_t & = 
  \begin{bmatrix}
   R(q_t)^T & 0 \\
   0 & R(q_t)^T
  \end{bmatrix}
  \begin{bmatrix}
   r^g \\
   r^b
  \end{bmatrix} \ud t + 
  \begin{bmatrix}
   \Sigma_W & 0 \\
   0 & \Sigma_W
  \end{bmatrix} \ud W_t,
\end{align*}
where the angular velocity is given by \cite{zamani2013thesis},
\begin{equation*}
 \omega_t = \Big( \sin(\frac{2\pi}{15}t),~-\sin(\frac{2\pi}{18}t+\frac{\pi}{20}),~\cos(\frac{2\pi}{17}t) \Big),
 \label{eq:angular_velocity}
\end{equation*}
$r^g=(0,0,-1)$ and $r^b=(1/\sqrt{2},0,1/\sqrt{2})$ represent the direction of the gravity 
and the local magnetic field, 
and $\Sigma_B$ and $\Sigma_W$ are $3\times3$ diagonal matrices 
associated with the process noise and the sensor noise, respectively.

The following filters are implemented for the comparison:
\begin{enumerate}
 \item IEKF: the algorithm is described in Sec. \ref{sec:IEKF}.
 \item MEKF: the algorithm is described in Sec. \ref{sec:MEKF}, using the modified Rodrigues parameter. 
 \item UKF: the algorithm is described in Sec. \ref{sec:UKF} and \cite{crassidis2003unscented}, 
            using the modified Rodrigues parameter.
 \item FPF-G: the FPF using the Galerkin gain functions, as described in Sec. \ref{sec:galerkin_gain} and Algorithm 2, 
              with fixed basis functions defined in Appendix \ref{sec:basis_fn}.
 \item FPF-K: the FPF using the kernel-based gain functions, as described in Sec. \ref{sec:kernel_gain} 
              and Algorithm 3, with $\epsilon=1$ and $K=10$.
\end{enumerate}

The performance metric is the root-mean-squared error (RMSE) 
\cite{zamani2013thesis, izadi2015comparison}:
\begin{equation*}
 \text{RMSE}_t = \sqrt{(1/M)\sum\nolimits_{j=1}^M \big( \delta\phi_t^j \big)^2},
 \label{eq:RMSE}
\end{equation*}
where $\delta\phi_t^j$ is the rotation angle error at time $t$ for the $j$-th Monte Carlo run, $j=1,2,...,M$. 
The rotation angle error is defined as follows: Let $q_t$ and $\hat{q}_t$ denote the true and estimated attitude at time $t$, and let 
$\delta q_t := \hat{q}_t^{-1}\otimes q_t$ represent the estimation error, 
then $\delta\phi_t = 2 \arccos(|\delta q_0|)$, where $\delta q_0$ is the first component of $\delta q_t$.
 
The filters are initialized with a ``concentrated Gaussian distribution'' \cite{wang2006error}, 
denoted as $N(q_I, \Sigma_0)$, whose mean $q_I$ is the identity quaternion, 
and $\Sigma_0$ is a diagonal matrix representing the variance in each axis of the Lie algebra. 
The particles in the FPF algorithms are sampled from this distribution as follows: 
First, one generates samples $\{v^i\}_{i=1}^N$ from 
the Gaussian distribution $N(0,\Sigma_0)$ in $\Re^3$. 
Then, the particles $\{R_0^i\}_{i=1}^N$ are obtained by $R_0^i = \exp (\sk{v^i})$, and converted 
to the quaternions $\{q_0^i\}_{i=1}^N$.

The simulations are conducted over a finite time-horizon $t\in[0,T]$ with fixed 
time step $\Delta t$. The process noise $\Sigma_B$ has standard deviation (std. dev.) of $5\,({}^{\circ}/\text{s})$. 
To avoid numerical instability due to large gain values, each of the first three measurement updates 
in all filters is implemented sequentially on a partition of the time step $\Dt$ with $N_f$ 
uniform sub-intervals. For FPF-G, $N_f=100$; For other filters, $N_f=20$ when $\Sigma_W$ is large, 
and $N_f=30$ when $\Sigma_W$ is small.
The relevant parameters are listed in Table 1. 

\begin{table}[h!] 
 \caption{\footnotesize{TABLE 1: SIMULATION PARAMETERS}}
 \centering
 {\normalsize
 \begin{tabular}{|c|c|c|c|c|}
  \hline
    $\Sigma_B$  &  $T$  &  $\Dt$  &  $N$  &  $M$  \\ \hline
    $0.008727^2\id$  &  2  &  0.01  &  200  &  100  \\
  \hline
 \end{tabular}
 }
 \label{table:simulation_params}
\end{table}

Fig. 1 illustrates the filter performance with different initial distribution and sensor noise. 
In Fig. 1 (a), $\Sigma_0 = 0.5236^2\id$, 
corresponding to the std. dev. of $\degg{30}$, and the target is initialized from the same distribution.  
The sensor noise $\Sigma_W = 0.01745^2\id$, i.e., the std. dev. is $\degg{10}$. 
In this case, all the filters have nearly identical performance. 

In Fig. 1(b), $\Sigma_0 = 1.0472^2\id$, corresponding to the std. dev. of $\degg{60}$, 
and the target is initialized with fixed attitude\,--\,rotation of $\degg{180}$ about the axis $(3,1,4)$. 
The sensor noise $\Sigma_W = 0.05236^2\id$, i.e., the std. dev. is $\degg{30}$.
These parameters indicate larger initial estimation error and uncertainty of the filters, and larger sensor noise. 
In this case, the FPF-K converges significantly faster than the other filters. 
% The UKF outperforms the EKF-based filters, but converges slower than the FPF.

When the initial estimation error is large, 
the Galerkin scheme yields significant error in computing the gain functions,  
and thus FPF-G is not included in Fig. 1(b). 
It is expected that one would require additional basis functions in this case. 
The computational complexity of the Galerkin scheme for one measurement update is approximately 
linear in the number of particles, whereas it is approximately quadratic for computing the kernel-based gain functions.

\section{Conclusion}
\label{sec:conclusion}

In this paper, the feedback particle filter was presented for the problem of attitude estimation.  
% a continuous-time nonlinear filtering problem on the Lie group $SO(3)$. 
The FPF is an intrinsic algorithm, possesses a gain-feedback structure and 
automatically respects the geometric constraint of the manifold. 
The algorithm was described using both the rotation matrix and the quaternion coordinate. 
The performance of FPF and its comparison with other attitude filters was illustrated 
by numerical simulations.

The continuing research includes improving the computational efficiency of the gain function approximation,  
% exploring structural connections between FPF and the EKF-based filters; 
and application of FPF for attitude estimation with more complicated models 
with e.g., translation and sensor bias.

\appendix
\label{sec:appendix}

\subsection{Basis functions in Galerkin scheme}
\label{sec:basis_fn}

For the Galerkin scheme presented in Sec. \ref{sec:galerkin_gain}, 
the following basis functions on $SO(3)$ are considered, expressed using the quaternion:
\begin{equation*}
 \begin{aligned}
 \psi_1(q) & = 2q_1q_0 , ~~\psi_2(q) = 2q_1q_0 , \\
 \psi_3(q) & = 2q_1q_0 , ~~\psi_4(q) = 2q_0^2-1. \\
 \end{aligned}
 \label{eq:basis_SO3_quat_mode1}
\end{equation*}
In order to compute the matrix $A$ and the vector $b$ in the Galerkin scheme, the formulae
for the action of $E_1,~E_2,~E_3$ on these basis functions are provided in Table 2.

\begin{table}[h!] 
 \caption{\footnotesize{TABLE 2: ACTION OF $E_n$ ON BASIS FUNCTIONS}}
 \centering
 {\normalsize
 \begin{tabular}{c|c|c|c}
  \hline
             & $E_1\cdot$       & $E_2\cdot$       & $E_3\cdot$      \\ \hline
    $\psi_1$ & $q_0^2-q_1^2$    & $-q_1q_2-q_3q_0$ & $-q_1q_3+q_2q_0$ \\ \hline
    $\psi_2$ & $-q_1q_2+q_3q_0$ & $q_0^2-q_2^2$    & $-q_2q_3-q_1q_0$ \\ \hline
    $\psi_3$ & $-q_1q_3-q_2q_0$ & $-q_2q_3+q_1q_0$ & $q_0^2-q_3^2$   \\ \hline
    $\psi_4$ & $-2q_1q_0$       & $-2q_2q_0$       & $-2q_3q_0$       \\
  \hline
 \end{tabular}
 }
 \label{table:lie_deriv_quat}
\end{table}

% \begin{remark}
%  Using the definition of unit quaternion \eqref{eq:quat_def}, 
%  the basis functions given in \eqref{eq:basis_SO3_quat_mode1} are equivalently written as, 
%  \begin{equation}
%   \begin{aligned}
%   \psi_1 & = \omega_1 \sin(\theta) , ~~\psi_2 = \omega_2 \sin(\theta) , \\
%   \psi_3 & = \omega_3 \sin(\theta) , ~~\psi_4 = \cos(\theta), \\
%   \end{aligned}
%   \label{eq:basis_SO3_theta_mode1}
%  \end{equation}
%  which are called basis functions of the {\em first mode}. It is straightforward to define 
%  basis functions of the {\em second mode}, 
%  \begin{equation}
%   \begin{aligned}
%   \psi_5 & = 2\psi_1 \psi_4 , ~~\psi_6 = 2\psi_2 \psi_4 , \\
%   \psi_7 & = 2\psi_3 \psi_4 , ~~\psi_8 = 2\psi_4^2-1. \\
%   \end{aligned}
%   \label{eq:basis_SO3_quat_mode2}
% \end{equation}
% The action of $E_1,~E_2,~E_3$ on the second-mode basis functions can be calculated 
% by the product rule and using Table-\Rom{1}.
% \end{remark}

\bibliographystyle{plain}
\bibliography{CDC_16}

\end{document}